\documentclass[english,draft]{article}
\usepackage{amssymb,amscd,amsthm,amsmath,babel,amsfonts}

\theoremstyle{plain}
\newtheorem{theorem}{Theorem}[section]
\newtheorem{lemma}[theorem]{Lemma}

\newtheorem{proposition}[theorem]{Proposition}

\theoremstyle{definition}

\newtheorem{definition}[theorem]{Definition}
\newtheorem{remark}[theorem]{Remark}

\newtheorem{example}[theorem]{Example}

\newtheorem{question}[theorem]{Question}
\newtheorem{problem}[theorem]{Problem}
\newtheorem{condition}[theorem]{Condition}
\numberwithin{equation}{section}

\newcommand\opn[2]{%
    \newcommand{#1}{\operatorname{#2}}}

\newcommand\NN{{\mathbb N}}
\newcommand\ZZ{{\mathbb Z}}

\newcommand\QQ{{\mathbb Q}}

\newcommand\frk{\mathfrak}

\newcommand\mm{{\frk m}}
\newcommand\pp{{\frk p}}

\opn\charac{char}
\opn\projdim{proj\,dim}
\opn\depth{depth}
\opn\rank{rank}
\opn\rlex{rlex}
\opn\LEX{Lex}
\opn\fine{end}
\opn\lcm{lcm}
\opn\modnuo{mod}
\opn\supp{supp}
\opn\Ass{Ass}
\opn\Proj{Proj}
\opn\Spec{Spec}
\opn\Soc{Soc}
\opn\reg{reg}
\opn\Md{Md}
\opn\dBorel{dBorel}
\opn\Borel{Borel}
\opn\Shad{Shad}
\opn\Span{span}
\opn\Hs{Hilb}
\opn\ini{in}
\opn\Gin{Gin}
\opn\grade{grade}
\opn\greg{greg}
\opn\Ht{ht}
\opn\GL{GL}

\opn\Ker{Ker}
\opn\coker{coker}
 
\opn\im{Im}
\opn\Hom{Hom}
\opn\Tor{Tor}
\opn\Ext{Ext}
\opn\id{id}
\opn\Id{Id}
\opn\Homst{^*Hom}
\opn\Extst{^*Ext}
\opn\Gamst{^*\Gamma}
\opn\Hst{^*H}

\opn\lk{lk}
\opn\st{st}

\let\:=\colon
\let\phi=\varphi

\newcommand\thh{^{th}}

\newcommand\pnt{{\raise0.5mm\hbox{\large\bf.}}}

\newcommand\ra{\rightarrow}

\newcommand\p{^\prime}
\newcommand\ppr{^{\prime\prime}}
\newcommand\ov{\overline}
\def\cocoa{\mbox{\rm
   C\kern-.13em o\kern-.07 em C\kern-.13em o\kern-.15em A}}
\newcommand\sat{^{sat}}
\opn\po{pol}
\opn\pol{^{\bf p}} 

\begin{document}
\noindent  
{\Large On non-standard graded algebras}

\vspace{1.5cm}
\noindent
{\bf Giorgio Dalzotto}. {\small  Universit\`a di Pisa, Dipartimento di
           Matematica, Largo Pontecorvo 5, 56127 Pisa, Italy.
           {\it e-mail}:  dalzotto@mail.dm.unipi.it} 

\noindent
{\bf Enrico Sbarra}. {\small Ruhr-Universit\"at Bochum, Fakult\"at f\"ur Mathematik
NA 3/32, 44780 Bochum, Germany. {\it e-mail}:
enrico.sbarra@ruhr-uni-bochum.de}

\vspace{1.5cm}

\noindent
{\small {\bf Abstract}

\noindent
Positively graded algebras are fairly natural objects which are arduous
to be studied. In this article we query quotients of non-standard
graded polynomial rings with combinatorial and commutative algebra methods.
}

\vspace{1cm}
{\scriptsize
\setcounter{tocdepth}{0}
\noindent CONTENTS
\setcounter{tocdepth}{9}

\contentsline {section}{\rm \numberline {1}Weighted generic initial
  ideals}{\rm 4}
\contentsline {subsection}{\numberline {1.1}Existence of the generic initial ideal}{5}
\contentsline {subsection}{\numberline {1.2}Weighted strongly stable ideals}{6}
\contentsline {section}{\rm \numberline {2}Prime Avoidance}{\rm 7}
\contentsline {section}{\rm \numberline {3}Regularity}{\rm 8}
\contentsline {section}{\rm \numberline {4}Lexicographic ideals}{\rm 11}
\contentsline {subsection}{\numberline {4.1}Hilbert functions}{11}
\contentsline {subsection}{\numberline {4.2}Lexifiable ideals}{12}
\contentsline {subsection}{\numberline {4.3}Polarization}{17}
}

\noindent
\section*{Introduction}
It has been shown in many works of the last three decades 
that combinatorial methods applied to Commutative Algebra and
Algebraic Geometry  are very effective. Most of these, though, deal with standard-graded polynomial
rings, i.e. with polynomial rings where the {\it weight} or degree
of all the variables involved is $1$.

\noindent
Other than the work \cite{A}, which is concerned with positively graded
algebras having a specific Hilbert function (see \cite{S}, Chapter 10
for a nice survey) and the introductory work \cite{BR}, there is a
small amount of literature about non-standard graded (or {\it
  weighted}) algebras, where the weight of a variable can be any positive integer.

\noindent
This lack persuaded us  to explore  the realm of weighted graded
algebras  with algebraic and combinatorial tools.

\noindent
Our work is divided in four sections, each centred on one topic and
inspired by some of the more significant results in the standard case: generic initial ideals,
Prime Avoidance, Castelnuovo-Mumford regularity and lexicographic
ideals.

The first section is dedicated to the definition and main properties
of generic initial ideals. We first study and describe the
automorphisms of a weighted polynomial ring, which are necessary for
the definition, the existence of generic initial ideals is discussed 
and a Borel-fixed type of property, i.e. fixedness under the action of a subgroup of the
group of the automorphisms is proved. Moreover, we define a combinatorial counterpart
of this property, i.e. being ``weighted strongly stable'', and prove that
generic initial ideals are enriched with it.

In the second section we recover an analogue of the homogeneous prime
avoidance Lemma (Lemma \ref{avoiding}), which grants the existence of
an almost-regular form of degree equal to the least common multiples of
the weights. This weaker statement is though enough to evince some
conclusions which generalize the known fact that depth does not change
after taking generic initial ideals with respect to the degree
reverse-lexicographic 
order (Proposition \ref{profo}). 
This is performed under the assumption that each weight is divisible
by the previous ones.

The third section takes into consideration the Castelnuovo-Mumford
regularity, which is defined in terms of  local cohomology. In the
first part we report some technical lemmata, which are useful in what
follows as computational tools. Then we prove that the regularity of
an ideal can be calculated using its graded Betti numbers, as it can
be done in the standard case but with a correction due to the weights and the number of
their occurrences.
This is achieved in Theorem
\ref{mainreg} and is essentially a paraphrase of what was performed in
\cite{Be} for pseudograded algebras. The section concludes with  a
result which predicts that Castelnuovo-Mumford regularity does not
change when taking generic initial ideal (Proposition
\ref{reggin}). Again, this is proved under the hypothesis that each weight
is divisible by the previous ones, and is false in general, as shown
in Example \ref{regmale}.

In the last section we deal with lexicographic ideals. In a standard
graded polynomial ring lexicographic ideals are enhanced with many
features and are well understood. As a consequence a certain amount of
information about an arbitrary homogeneous ideal can be gathered by
studying the associated lexicographic ideal. 
In other contexts though, the generalization
of this notion turned out to be rather complex (cf. \cite{ADK}, 
\cite{DH} and \cite{MP}). In our situation there is a natural way to
define lexicographic ideals, but it is difficult to describe them and
to give a criterion to decide whether a homogeneous ideal in a given
weighted polynomial ring is {\it lexifiable}, i.e. admits an
associated lexicographic ideal with the same Hilbert function.

\noindent
First we recollect some known results about Hilbert functions of a
positively graded algebra and underline which are the factors which
make our analysis difficult by means of some examples. In particular
the {\it shadow} of a lexsegment is a lexsegment in the standard
graded case, but this fact does not hold in general in the
non-standard setting. We thus proceed
by proving Proposition \ref{chebene}, which yields a method to verify
whether a given  ideal is lexicographic. It is indeed enough to check
if  finitely many graded components are generated by lexsegments. 
This is accomplished by means of an invariant $G(w)$, which was
introduced in \cite{D}. From a computational point of view the test is
not optimal yet, since $G(w)$ not only depends on the weights but
grows rapidly with the number of variables having the same weight.

\noindent
 Next, we are interested in describing the Hilbert function of
lexicographic ideals, but to give a complete solution (cf. Problem 
\ref{problema}) is an hard task. Theorem \ref{pesiprimi} provides an
exhaustive answer for polynomial rings in two variables. The next
topic we handle is expressed by Question \ref{acaulay}. One would like
to know which sets of weights make a polynomial ring {\it
  Macaulay-Lex}, i.e. such all of its ideals are lexifiable. Theorem
\ref{pesi1} and the subsequent examples provide partial answers to
this issue. We would like to observe that Theorem \ref{pesiprimi} and
Theorem \ref{pesi1} provide a complete description of lexicographic
ideals in two variables. Still, even in two variables it is not clear
which ideals are lexifiable.

\noindent
As a final remark, which completes and concludes this survey, we also
mention the technique of {\it polarization}. In this setting
admissible numerical functions need not to be Hilbert functions of 
lexicographic ideals. Thus completely polarized ideals (which in the standard
case characterize lexicographic ideals), might be the right tool for a
Theorem {\em \`a la} Macaulay \cite{M}.

\medskip
\noindent
The calculations underlying many of the examples and the material of
the last section were carried out
using \cite{CoCoA}. We implemented some procedures (for the computation of 
Hilbert functions, generic initial ideals, polarization and
associated lexicographic ideals) which can be obtained by any of the
two authors. 

\section*{Notation}
In this paper we use some non-standard notation (!) which we illustrate
here.
When we consider polynomial rings with a non-standard grading, 
we mean that we work over an infinite field $K$ of characteristic $0$
and assume  the degrees of the variables to be positive integers with no further
restriction. 

\noindent
We order the variables by increasing degree or {\it weight} and often group together those
with the same degree. Therefore we denote the
polynomial ring by  $R=K[{\bf  X_1},\ldots ,{\bf X_n}]$, where 
${ \bf X_i}=(X_{i1}, \ldots , X_{il_i })$, $\deg X_{ij}=q_i$
for any $j=1,\ldots,l_i$, and  $q_1< q_2<\ldots< q_n$.

\noindent
It is convenient to denote by $R_{[i]}$ the polynomial
ring $K[{\bf  X_1},\ldots ,{\bf X_i}]$.
We let $w$ be the weight vector $(\deg X_{11}, \ldots, \deg
X_{nl_n})$ so that  $(R,w)$ stands for a polynomial ring with the
graduation given by $w$. If $w$ does not play an explicit
role we denote $(R,w)$ simply by $R$. 

\noindent If not elsewhere specified we consider term orderings $>$ which
are degree compatible and assume $X_{ij}>X_{ik}$ if $j<k$, $i=1,\ldots,n$.

\noindent
Since they are often used, it
may be convenient to fix some notation for the total numbers of
variables and the least common multiple of the weights. Thus we let
$l=\sum_{i=1}^n l_i$ and $q=\lcm(q_1,\ldots,q_n)$. 

\noindent
Finally, given a set $A\subseteq R_d$, 
$\langle A \rangle$ denotes the $K$-vector space spanned by $A$.
If $V\subseteq R_d$ is a  $K$-vector space, $\{V\}$ denotes its monomial basis.

\section{Weighted generic initial ideals}\label{wgin}
Generic initial ideals are monomial ideals attached to  homogeneous 
ideals. It has been shown in many works that generic 
initial ideals, although simpler in nature, still contain a
considerable amount of information about the original geometrical
object. In order to extend their definition in our setting we have 
first of all to  understand which are the graded automorphisms of $R$. 

\begin{proposition} 
The assignment
$$\phi(X_{ij})=\sum_{h=1}^{l_i}a^i_{jh}X_{ih}+\psi_{ij}({\bf
  X_1},\ldots,{\bf X_{i-1}}),$$  
where, for all $i,j$,  $\psi_{ij}$ are homogeneous polynomials in $R_{[i-1]}$ of degree
$q_i$ and $A_i \doteq (a^i_{jh})_{j,h=1,\ldots,l_i}\in M_{l_i}(K)$ are
invertible matrices, defines a graded automorphism of $R$. Vice versa,
any graded automorphism of $R$ is of this kind.
\end{proposition}
\begin{proof}
Since $q_i<q_j$ if $i<j$, the requirement that $\phi$ is a graded
homomorphism forces $\psi_{ij}$ to be polynomials in the first $i-1$
sets of variables. Thus it is sufficient to show that $\phi$ is
surjective if and only if $A_i$ are invertible for all $i=1,\ldots,n$.
With  some abuse of notation we write $\phi({\bf X_i})=A_i{\bf
  X_i}+B_{i-1}$, where $B_{i-1}$ is a $l_i\times 1$ matrix with 
entries in $R_{[i-1]}$.
If $\phi$ is surjective then for all $i=1,\ldots,n$ there exists a
$l_i\times 1$ matrix $C_i$ with entries in $R_{[i]}$ such 
that ${\bf X_i}=\phi(C_i)$. If we write $C_i$ as $D_i{\bf X_i}+E_{i-1}$,
where $D_i \in M_{l_i}(K)$ and $E_{i-1}$ is a $l_i\times 1$ matrix with
entries in $R_{[i-1]}$, we get that ${\bf X_i}=D_iA_i{\bf X_i}+F_{i-1}$
where $F_{i-1}$ consists of polynomials in ${\bf X_1},\ldots,{\bf
  X_{i-1}}$.
Therefore $F_{i-1}=0$ and $D_iA_i=I$. Vice versa, suppose that $A_i$ is invertible
for $i=1,\ldots,n$. Since $A_1$ is invertible, $\phi(A_1^{-1}{\bf
  X_1})={\bf X_1}$. If $i>1$, we have $\phi(A_i^{-1}{\bf X_i})={\bf
  X_i}+A_i^{-1}B_{i-1}$, where $A_i^{-1}B_{i-1}$ has entries in
$R_{[i-1]}$. 
 Thus there exists a $l_i\times 1$ matrix $C_{i-1}$ with entries in 
 $R_{[i-1]}$ such that
$\phi(C_{i-1})=A_{i}^{-1}B_{i-1}$ and ${\bf X_i}=\phi(A_i^{-1}{\bf
  X_i})-\phi(C_{i-1})=\phi(A_i^{-1}{\bf X_i}-C_{i-1})$, as required.
\end{proof}

Let  $T$ be the subset of {\it upper triangular automorphisms} 
consisting of those graded automorphisms of $R$ such that, for all $i=1,\ldots,n$, $A_i$
is an upper triangular invertible matrix. By the previous proposition 
it is clear that $T$ is a group, since  upper triangular invertible
matrices form a group (which is called {\it Borel group}).
Let  $U$ be the set of
the {\it elementary upper triangular automorphisms} $\tau_{ij}^{rc}$,
where $r<j$ and $c\in K$, determined by the assignment  $\tau_{ij}^{rc}(X_{ij})=X_{ij} +
cX_{ir}$ and $\tau_{ij}^{rc}(X_{hk})=X_{hk}$ if $(h,k)\neq (i,j)$.
Finally, let $N$ be the set of {\it  elementary non-linear automorphisms}
$\eta_{ij}^m$, where $m$ is a term of degree $q_i$ in $R_{[i-1]}$
defined by 
$\eta_{ij}^m(X_{ij})=X_{ij}+m$ and $\eta_{ij}^m(X_{hk})=X_{hk}$ if
 $(h,k)\neq (i,j)$. 

\begin{proposition}\label{generato}
$T$ is generated by the diagonal subgroup, by $U$ and by $N$.
\end{proposition}
\begin{proof}
The proof is an easy induction on the number of variables $l=l_1+\ldots + l_n$.
If $l=1$, the only graded automorphisms are the diagonal automorphisms.
Let $\phi \in T$, we say
$\phi(X_{ij})=\sum_{h=1}^{l_i}a^i_{jh}X_{ih}+\psi_{ij}$, and
$A_i\doteq (a_{jh}^i)_{j,h=1,\ldots,l_i}$ are upper triangular
  invertible matrices;
also, let
 $\phi\p$ be 
defined by  $\phi\p(X_{nl_n})= X_{nl_n}$ and $\phi\p(X_{ij})=
\phi(X_{ij})$ if $(i,j) \neq (n,l_n)$. It is clear that $\phi\p$ is an
automorphism and that therefore belongs to $T$. We now write the polynomial
$\psi_{nl_n}$ as sum of monomials $m^h_{nl_n}\in K[{\bf X_1},\dots ,
{\bf X_{n-1}}]$, $h=1,\ldots,s$, of degree $q_n$.
 We want to find a decomposition of
$\phi$ by means of elementary non-linear, upper triangular, diagonal
automorphisms and of $\phi\p$. This leads to the conclusion by
induction, since $\phi\p$ fixes the last variable and can be thought
of as an automorphism of a polynomial ring in $l-1$ variables. We denote
by $\delta_{ij}^c$, with $c\in K$, the diagonal automorphism
defined by $\delta^c_{ij}(X_{ij})=cX_{ij}$ and
$\delta_{ij}^c(X_{hk})=X_{hk}$ if $(h,k)\neq (i,j)$. Moreover,  since
$a^n_{l_nl_n}\neq 0$, we may write 
$u_h\doteq m_{nl_n}^h/a^n_{l_nl_n}$ for all $h=1,\ldots,s$ and  $b_k\doteq
\frac{a_{nk}^n}{a_{l_nl_n}^n}$, for $k=1,\ldots,l_n-1$. It is now
easy to see that 

$$\phi=\eta_{nl_n}^{u_1} \circ
\cdots \circ  \eta_{nl_n}^{u_s} \circ \tau_{nl_n}^
{1b_1} \circ \cdots \circ \tau_{nl_n}^ 
{l_n -1b_{l_n-1}} \circ
\delta_{nl_n}^{a_{l_nl_n}^n} \circ  \phi\p,$$
as desired.
\end{proof}

\subsection{Existence of the  generic initial ideal} In the standard graded
case the generic initial ideal $\Gin(I)$ of an ideal $I\subset
K[X_1,\ldots,X_n]$
plays a central role in problems
regarding Hilbert functions and free resolutions of graded
ideals. Since $\Gin(I)$ with respect to some assigned term order is
defined as the initial ideal of $gI$, where $g$ is a generic change of
coordinates, i.e. that $g$ is a matrix chosen out of a
Zariski non-empty open set of $\GL_n(K)$, one way of computing it
is the following. Let the $n^2$ entries of $g$ be new
indeterminates, we say $Y_{ij}$, with $i,j=1,\ldots,n$. Write $gI$
explicitly and apply the Buchsberger's Algorithm to compute $\ini(I)$ as
an ideal of 
$K(Y_{ij})[X_1,\ldots,X_n]$. After finitely many computations of the so-called
{\it $S$-pairs} the process finishes, the output result is the sought
after  monomial ideal - in the variables
$X_i$ only - and the Zariski open set consists of all those matrices for which
the finitely many polynomial denominators of the $S$-pairs are
non-zero. If one considers this point of view, it is evident that
weights do not play any role in the construction, which is thus also
possible in the weighted case. Thus we can talk of generic initial
ideals of homogeneous ideals in a non-standard graded algebra.

\noindent In the standard case it is well-known that generic
initial ideals are {\it Borel-fixed}, i.e. fixed under the action of
the Borel subgroup of $\GL_n(K)$ consisting of the upper triangular
invertible matrices. 

\begin{theorem}\label{fissita}
Let $I$ be a homogeneous ideal of a weighted polynomial algebra $R$.
Then $\Gin(I)$ is  $T$-fixed, i.e.  $\phi(\Gin(I))=\Gin(I)$ for all
$\phi \in T$.
\end{theorem}
\begin{proof}
We only need to observe that by applying a non-linear automorphism to
a monomial 
$u$ one obtains a polynomial of the form $u+v$ where $v$ is bigger
than $u$ in the chosen term order. By Proposition \ref{generato}, this
is enough to argue as in the standard case, see for instance the proof
of Theorem 15.20 in \cite{E}.
\end{proof}

\subsection{Weighted strongly stable ideals} 
Generic initial ideals in a standard graded polynomial ring are
characterized combinatorially, the simplicity of this description
depending on the characteristic of the base field. In a weighted
polynomial ring over a base field of characteristic $0$ the same can
be performed, via the following definition.

\begin{definition}\label{ss}
Let $I$ be a monomial ideal. $I$ is called {\it (strongly) stable} if
the following holds:  for every $u\in I$, if $X_{ij}\mid u$ then  $\frac{X_{ih}u}{X_{ij}}\in I$, 
for every  $h<j$ and $\frac{vu}{X_{ij}}\in I$ for all monomials $v$ of
degree $q_i$ in $R_{[i-1]}$. 
\end{definition}

It is not difficult to prove that weighted generic initial ideals are
stable according to this definition.

\begin{proposition}\label{stabilita}
Let $I$ be an homogeneous ideal. $I$ is $T$-fixed if and only if $I$
is strongly stable. 
\end{proposition}
\begin{proof}
One begins by observing that $I$ is fixed by the subgroup of diagonal
matrices if and only if $I$ is monomial. Let $m=X_{ij}^tm\p$, where
$X_{ij}\nmid m\p$. The images $\tau_{ij}^{rc}(m)$, $\eta_{ij}^s(m)$,
with $\deg s=q_i$ can be written as $(X_{ij}+cX_{ir})^tm\p$ and 
$(X_{ij}+s)^tm\p$ respectively. If $I$ is $T$-fixed both polynomials,
and so each of their monomials, belong to $I$. In particular the
conditions which define strongly stable ideals are
verified. Conversely, if $I$ is strongly stable the same argument shows
that $I$ is fixed by the action of the generators of $U$ and $N$, and
is therefore $T$-fixed.
\end{proof}

One of the key properties of strongly stable ideals in the standard
graded polynomial ring $K[X_1,\ldots,X_n]$ is that
$I:(X_1,\ldots,X_n)=I:X_n$. In fact, beside the trivial inclusion
$I:(X_1,\ldots,X_n)\subseteq I:X_n$ one has that $m\in I:x_n$ iff
$mX_n\in I$, which, because of the stability property, implies $mX_i\in
I$ for all $i$. It is quite clear how this property is weakened in the
more general case where variables might have different weights. In
particular the good property of stable ideals with respect to taking
colons with the last variables plays a central role in the
construction of the Eliahou-Kervaire resolution \cite{EK} of such an
ideal. This is a completely described  minimal graded free resolution 
of such an ideal in terms of its minimal set of monomial generators.
On the other hand being able to construct such a resolution having no 
restriction on the weight vector would mean to know how to describe a
minimal resolution of any monomial ideal, since given any
such ideal $I$, one can choose weights so that in the corresponding
polynomial ring $I$ is stable, as the next example shows.

\begin{example}
Let $A$ be a set of monomials in $n$ variables. Then, there exist
non-negative integers $q_1,\ldots,q_n$ such that in the
weighted polynomial 
ring $(K[X_1,\ldots,X_n],(q_1,\ldots,q_n))$  the ideal generated by
$A$ is strongly stable. In fact, one can choose 
weights in such a way that none of the  exchanges which were described in Definition
\ref{ss} is possible. For instance, it is enough to choose
$q_1<q_2<\ldots<q_n$ so that $2q_1>q_n$, we say $q_1=n+1,q_2=n+2, \ldots,q_n=2n$.
\end{example}

\section{Prime Avoidance}\label{primea}
A simple fact of linear algebra gives rise to a powerful tool when
combined with techniques dealing with generic forms. This is known as
Homogeneous 
Prime Avoidance: If
$\pp_1,\ldots,\pp_n$ are prime ideals strictly contained in the
graded maximal ideal of a standard graded algebra over an infinite
field then there exists a homogeneous form of
degree $1$ in $\mm\setminus \cup_i\pp_i$. It turns out to be essential
in many proofs, since avoiding a finite number of primes is an open
property.

\begin{lemma}[Weighted Prime Avoidance]\label{avoiding}
Let   $q=\lcm (q_1,\ldots,q_n)$ and let $\pp_1,\ldots,\pp_n$ be prime ideals with $\pp_i
\subsetneq \mm$. Then
$\mm_q\setminus \cup_i (\pp_i)_q \neq \emptyset$.
\end{lemma}
\begin{proof}
Since the prime ideals are strictly contained in the maximal ideal, we
have that $(\pp_i)_q \neq \mm_q$ for all $i$. Else, one would have
that  $(\pp_i)_q
= \mm_q$ and  $X_{jk}\in \pp_i$, for all $j$ and $k$, since $X_{jk}^{q/q_j}\in \mm_q$ and
$\pp$ is prime. But the infinite vector space $\mm_q$ cannot be
written as a finite union of proper subspaces $\pp_q$, and the claim follows.
\end{proof}

The next example shows that in general it is not possible to find such
a form in a smaller degree.

\begin{example}
Let $(R,w)=(K[X,Y],(2,3))$. If $\pp_1=(X)$ and
$\pp_2=(Y)$ then the smallest degree $d$  such that $(X,Y)_d\supsetneq (\pp_1)_d\cup(\pp_2)_d$ is $6$. 
\end{example}

In the following we recover some results which are known
in the standard case, provided that some condition on the
weights is assumed. It may be convenient to state one of these conditions here.

\begin{condition}\label{multipli}
$(R,w)$ is a weighted polynomial ring with $q_i\mid q_{i+1}$ for
$i=1,\ldots,n-1$.
\end{condition}

\begin{lemma}
Let $(R,w)$ be a ring for which Condition \ref{multipli} is satisfied,
and let $I\subseteq R$ be a strongly stable ideal. For any $i=1,\ldots , n$
and $j=1,\dots,l_i$, one has $$I:X_{ij}^\infty = I:(X_{11}, \ldots, X_{ij})^\infty .$$
\end{lemma}

\begin{proof}
We only have to prove the inclusion $\subseteq$ since the other one
is obvious. Let $m$ be a monomial such that $mX_{ij}^s \in I$ for some
$s\in\NN$. 
Since $I$ is strongly stable, $mX_{ih}^s \in I$ for
any $1\leq h \leq j$; furthermore the assumption on the degrees of
the indeterminates implies that  $\deg X_{ij} = q_i = q_h \frac{q_i}{q_h}= \deg
X_{hk}^{q_i/ q_h}$, and consequently $m (X_{hk}^{q_i/q_h})^s \in I$
for any $1\leq h\leq i-1$ and $1\leq k\leq l_h$, as desired.
\end{proof}

As a consequence we obtain the following proposition.

\begin{proposition}\label{ultimoindice}
Let $(R,w)$ be a ring for which Condition \ref{multipli} is satisfied,
let 
$I$ be a strongly stable ideal and let $X_{hk}$ be the (lex-)smallest
variable which divides some minimal generator of $I$. Then
$X_{hk+1},\ldots,X_{nl_n}$ form a maximal regular sequence on $R/I$.
\end{proposition}
\begin{proof}
Clearly the elements $X_{hk+1},\ldots,X_{nl_n}$ form a regular sequence on $R/I$. 
Since in the quotient ring $\ov{R}=K[X_{11}, \dots,
X_{hk}]$ the ideal $\ov{I}$ is strongly stable, by the previous
lemma $\overline{I}\sat= \overline{I}:X_{hk}^\infty \neq
\overline{I}$, which implies that $\depth \ov{R}/\ov{I}=0$.
\end{proof}

We recall the following theorem \cite{BS}, which holds independently
of the given weights and is needed for the proof of  the final result of this section.

\begin{theorem}\label{iein}
Let $F$ be a free $R-$module with basis and consider the degree
reverse lexicographic 
monomial order. Let $M$ be a graded submodule of $F$. The elements
$X_{nl_n},X_{nl_n-1}, \ldots,X_{ij+1}, X_{ij}$ form a regular sequence on $F/M$ if and
only if they form a regular sequence on $F/\ini(M)$.
\end{theorem}
\begin{proof}
See that of Theorem 15.13 in \cite{E}.
\end{proof}

\begin{theorem}\label{profo}
Let $(R,w)$ be a ring for which Condition \ref{multipli} is satisfied,
and consider the degree reverse lexicographic order. Then, for any homogeneous ideal  $I\subseteq R$,
$$\depth R/I=\depth R/\Gin(I).$$
\end{theorem}
\begin{proof}
Since $\depth R/I \geq \depth R/\Gin(I)$, we may assume $\depth
 R/I>0$. By Lemma \ref{avoiding}, a generic form of degree $q_n$ is
 a non-zerodivisor on $R/I$. Thus, after a generic change of coordinates, we may assume that
$X_{nl_n}, X_{nl_n-1},\ldots, X_{ij}$ is a maximal $R/I$-regular
 sequence and
$\Gin(I)=\ini(I)$. By Theorem \ref{fissita} and Proposition \ref{stabilita}
$\Gin(I)$ is strongly stable, and consequently, by Proposition
 \ref{ultimoindice}, 
there is a maximal $R/\Gin(I)$-regular sequence $X_{nl_n}, X_{nl_n-1},
\ldots, X_{hk}$. Now Theorem \ref{iein} yields that $(h,k)=(i,j)$,
 from which the conclusion is straightforward. 
\end{proof}


\section{Regularity}\label{regu}
Local cohomology modules of a graded module over a weighted polynomial ring have a graded structure
arising from resolutions by graded injective modules or equivalently
from the construction of the 
\v{C}ech complex. The usual definition of Castelnuovo-Mumford
regularity by means of local cohomology 
still works in this context and we recall it here.
Let $H_\mm^i(M)$ denote the $i\thh$ graded local cohomology module of
the graded $R$-module $M$ with  support on the graded maximal ideal $\mm$. 

\begin{definition}
Let $R$ be a weighted polynomial ring with graded maximal ideal $\mm$. We let 
$$a^i(M)\doteq \left\{\begin{array}{ll}\max\{j\in\ZZ\: H_\mm^i(M)_j\neq 0\}& \hbox{if } H^i_\mm(M)\neq 0\\ -\infty & \hbox{otherwise}\\
                      \end{array}
\right .$$ 
denote the end of the $i\thh$ local cohomology module of $M$. The {\em Castelnuovo-Mumford regularity} of $M$ is then
$\reg M =\max_{1\leq i\leq \dim M}\{a^i(M)+i\}$.
\end{definition}

However, one of the aspects that made the Castelnuovo-Mumford
regularity interesting, i.e. its direct interpretation through the Betti
numbers of the minimal free resolution by means of the formula
\begin{equation}\label{zuca}
\reg M =\max_{i\geq 0}\{b_i(M)-i\},
\end{equation}
where $b_i(M)\doteq\max_{j\in\ZZ}\{\beta_{ij}(M)\neq 0\}$, fails in the general weighted
case. In this section we re-prove some results about regularity which
still hold in the weighted case, in order to give in Theorem
\ref{mainreg} a formula that
generalizes \eqref{zuca}. In the last
part, we consider the regularity of a generic initial ideal $\Gin(I)$
and prove in Proposition \ref{reggin}
that under some assumption on  weights it does not differ
from that of $I$. Also, we provide a counterexample that shows that in
general there is no analogue of the well-known theorem \cite{BS}
valid in the standard case.


We start by recalling some lemmata which are useful in order to
control regularity in the non-standard case.

\begin{lemma}\label{shortexact}
Let $0 \ra N \ra M \ra Q \ra 0$ be a short exact sequence of finitely
generated graded $R$-modules. Then 
\begin{itemize}
\item[(i)]
$\reg N  \leq \max\{\reg M, \reg Q+1\}$;
\item[(ii)]
$\reg M  \leq \max\{\reg N, \reg Q\}$;
\item[(iii)]
$\reg Q \leq \max\{\reg N-1, \reg M\}$;
\item[(iv)]
If $N$ has finite length, then $\reg M =\max\{\reg N, \reg Q\}$.
\end{itemize}
\end{lemma}
\begin{proof}
The proofs of $(i)-(iii)$ are easy and  descend from the
use of the long exact sequence in cohomology $\ldots \ra
H_\mm^{i-1}(Q) \ra H_\mm^i(N) \ra H_\mm^i(M) \ra \ldots$. 

As for the proof of $(iv)$, it is clear that $\reg N=a^0(N)$
and $a^0(M)=\max\{a^0(N),a^0(Q)\}$. Thus,
\begin{equation*}\begin{split}
\reg M \doteq&\max\{a^0(M),\max_{i>0}\{a^i(M)+i\}\}\\
=&\max\{a^0(N), a^0(Q),\max_{i>0}\{a^i(Q)+i\}\},
\end{split}
\end{equation*}
 as desired.
\end{proof}

\begin{lemma}\label{modulo}
Let $M$ be a finitely generated graded $R$-module and let $x\in
R_d$. If $x$ is a non-zerodivisor on $M$
then $\reg M/xM=\reg M +(d-1)$. More generally, if $x$ is such
that $(0:_Mx)$ has finite length, then
$$\reg M =\max\{\reg 0:_Mx,\reg M/xM -(d-1)\}.$$
\end{lemma}
\begin{proof}
From the exact sequence $0\ra (0:_Mx)(-d) \ra M(-d) \ra M \ra M/xM \ra
0$ one obtains the two short exact sequences $0\ra (0:_Mx)(-d) \ra
M(-d) \ra xM \ra 0$ and $0 \ra xM \ra M \ra M/xM \ra 0$ so that the
proof follows easily as an application of  Lemma \ref{shortexact}.
\end{proof}

\begin{lemma}
Let $x\in R_d$ such that $0:_Mx$ is of finite length. Then for all
$i\geq 0$
$$a^{i+1}_\mm(M)+d\leq a^i_\mm(M/xM)\leq \max\{a^i_\mm(M), a^{i+1}_\mm(M)+d\}.$$
\end{lemma}
\begin{proof}
From the two short exact sequences contained in the proof of the last lemma
we deduce that $H^i_\mm(M(-d))\simeq H^i_\mm(xM)$ for all $i>0$ and
obtain the long exact sequence in cohomology
$\ldots\ra H^i_\mm(M) \ra H^i_\mm(M/xM) \ra H^{i+1}_\mm(xM) \ra
H^{i+1}_\mm(M) \ra\ldots$.
To prove the first inequality, it is enough to observe that, if
$a^i_\mm(M/xM)< a^{i+1}_\mm(M)+d$ the above long exact sequence in
degree $a^{i+1}_\mm(M)+d$ would deliver a contradiction. The proof of
the second inequality is analogous.
\end{proof}

Let $M$ be a finitely generated $R$-module of finite projective
dimension $s$.
For $i=1,\ldots,s$ let as before $b_i(M)\doteq\max_{j\in\ZZ}\{\beta_{ij}(M)\neq 0\}$.

\begin{theorem}\label{mainreg}
Let $R=K[{\bf  X_1},\ldots ,{\bf X_n}]$ be a graded polynomial ring and
let $M$ be a finitely generated $R$-module with $\projdim M<\infty$.  Then
$$\reg M =\max_{i\geq 0}\{a^i_\mm(M)+i\}=\max_{i\geq 0}\{b_i(M)-i\}-\sum_{j=1}^nl_j(q_j-1).$$
\end{theorem}
\begin{proof}
By virtue of the previous lemma and induction on the number of
variables one first proves that 
\begin{equation}\label{induco}
b_0(M)\leq \max_{i\geq
  0}\{a^i_\mm(M)+i\}+\sum_{j=1}^nl_j(q_j-1).
\end{equation}
Moreover, it is easy to verify that, if $F$ is a free $R$-module, 
\begin{equation}\label{libero}
a_\mm^n(F)=b_0(F)-\sum_{i=1}^nl_iq_i.
\end{equation}
The  assertion follows by the use of \eqref{induco} and \eqref{libero}
combined with an induction argument on the projective dimension of
$M$. The proof is an adaptation of that of  Theorem 5.5 in \cite{Be} 
(to which the interested reader is referred) and, therefore, the
details are omitted here.
\end{proof}

\begin{proposition}\label{reggin}
Let $(R,w)$ be a weighted polynomial ring for which Condition
\ref{multipli} is satisfied and consider the  degree reverse 
lexicographic order. If $I$ an homogeneous ideal of
$R$ then $$\reg R/I=\reg R/\Gin(I).$$
\end{proposition}
\begin{proof}
Since $q_n=\lcm(q_1,\ldots,q_n)$, by Lemma \ref{avoiding} a
generic form  in $\mm_{q_n}$ does not belong to any associated prime
$\pp\neq\mm$ of $I$. By applying a generic automorphism we may assume that
$\Gin(I)=\ini(I)$ and that $X_{nl_n}$ is almost-regular. Therefore
$(I:X_{nl_n})/I$ and, consequently, $\ini(I:X_{nl_n})/\ini(I)\simeq
(\ini(I):X_{nl_n})/\ini(I)$ have finite length. By Lemma \ref{modulo},
it is enough to verify that
$\reg (I:X_{nl_n})/I=\reg (\ini(I):X_{nl_n})/\ini(I)$ in order to
apply induction on the numbers of the variables, since the assumption
on the weights still holds for  $R/(X_{nl_n})$. But this is clear
because the above modules coincide with their $0\thh$ local cohomology
module and they have the same Hilbert function.
\end{proof}

Notice that the assumption on the weights is essential in order to
have a generic form of the right degree for the induction. The
following example shows that the above result cannot be extended for
any choice of  weights.

\begin{example}\label{regmale}
Let $(R,w)=(K[X,Y,Z],(2,4,5))$ and $I=(XY,YZ,X^5)$ and consider the
degree reverse lexicographic order.  A
computation with \cite{CoCoA} shows that $\Gin(I)=(X^3, X^2Z, XY^2,
Y^3Z)$. By Theorem \ref{mainreg}, the regularity of $I$ and $\Gin(I)$
can be computed by the use of the resolutions
$$0 \ra R(-14)\oplus R(-11) \ra R(-10)\oplus R(-9)\oplus R(-6) \ra  I
\ra 0$$ and 
\begin{equation*}
\begin{split}
0 \ra  R(-19) &\ra  R(-19)\oplus R(-17)\oplus R(-14)\oplus R(-11) \ra\\ 
& \ra R(-17)\oplus R(-10)\oplus R(-9)\oplus R(-6) \ra  \Gin(I) \ra 0\\
\end{split}
\end{equation*}
of $I$ and $\Gin(I)$ respectively.

This example points out that also Proposition \ref{profo} is not valid
without Condition \ref{multipli}.
\end{example}

\section{Lexicographic ideals}\label{lexico}
Although the definition and some of the main properties of Hilbert
functions are still valid in a non-standard setting, a great deal is
still unknown about them. In particular Macaulay's Theorem, which
provides a necessary and sufficient condition for a numerical function
to be the Hilbert function of a finitely generated standard graded
algebra has no counterpart in the weighted case. The main tool which
is involved in this context, lexicographic ideals, can be easily
defined in the non-standard case, but they are not so easily
investigated, as the following analysis shows.

\subsection{Hilbert functions}
Here we point out some facts about non-standard graded algebras 
which are relevant for our purposes. We start by recalling the 
well-known Hilbert-Serre Theorem: Let $I$ be a homogeneous ideal in
$(R,w)$. The Poincare series ${\it P}(R/I,t)$ of $R/I$ is a rational function in
$t$ of the form $g(t)/\prod_{i=1}^n (1-t^{q_i})^{l_i}$, where $g(t) \in
\ZZ[t]$.
\noindent
It is known that  the Hilbert function of $R/I$ is 
quasi-polynomial. Some more information is provided by the following
result to be found in \cite{B}, Theorem 2.2.

\begin{proposition}
Let $I$ be a homogeneous ideal in $(R,w)$ and let $d$ be the order of
the pole of ${\it P}(R/I,t)$ at the point $t=1$. Then there exist $q=\lcm(q_1,\ldots,q_n)$
polynomials $p_0,\ldots,p_{q-1}\in \QQ[t]$ of degree at most $d-1$ with coefficients in
$[q^{d-1}(d-1)!]^{-1}\ZZ$ such that, for all $l \gg 0$,
$$H_{R/I}(l)=p_j (l) \hbox{    for } l\equiv j\modnuo q$$
\end{proposition}

It is also worth observing that in general some of the Hilbert
polynomials described in the above proposition can be $0$. Also in the
case $\gcd(q_1,\ldots,q_n)=1$ the vanishing of the Hilbert function of
$R/I$ in $t=t_0$ does not imply that ${\it H}_{R/I}(t)=0$ for all
$t>t_0$. However, this is true for $H_R(t)$ if $t_0$ is bigger than the Frobenius
number of $q_1,\ldots,q_n$ (cf. \cite{SS}, Chapter 1, Section 3 for
more details about this subject).

\begin{remark}\label{tuttouguale}
Let $(R,w)$ be a weighted polynomial ring. If $w_i=q$ for all $i$,
then the Hilbert function $H_{(R,w)}(t)$ is equal to
$H_{(R,(1,\ldots,1))}(t/q)$ if $q\mid t$ and $0$
otherwise; this case is thus essentially equivalent to the standard
case. The same observation shows that one may assume without loss of
generality that the $\gcd$ of the weights is $1$.
\end{remark}


\noindent
Another pathology of the weighted case is shown in the following
example \cite{BR}.

\begin{example}
Let $(R,w)=(K[X,Y,Z,T],(1,6,10,15))$.

\noindent
The monomial $XY^4Z^2T$ has degree $60$, 
but it is not multiple of any monomial of degree $30$.
\end{example}

However, one can show that this can only occur in low degrees, as it
is shown in \cite{BR}, Proposition 4B.5. One makes use of an
invariant introduced in \cite{D}, which we denote by  $G(w)$. For the reader's sake
we recall here the result, omitting the definition of $G(w)$ since
it is not essential in what follows. 

\begin{proposition}\label{beltrametti}
Let $(R,w)$ be a weighted polynomial ring and let $n>G(w)$. 
Then every monomial of $R_{n+hq}$ is divisible by a
monomial in $R_{hq}$, for any $h\in \NN$.
\end{proposition}

One might wonder if the same holds for an arbitrary ideal generated in more
than one degree, i.e. if  there exists $l\in \NN$ such that
for all $r\gg 0$ one has $I_r=I_lR_{r-l}$. Unfortunately this is false
and it partially explains why the study of lexicographic ideals is complicated.

\begin{example}
Let $(R,w)=(K[X,Y,Z],(2,2,3))$ and $I=(X^\alpha,XYZ)$ for some
integer $\alpha>1$. Let us suppose that there exists some $l\in
\NN$ such that $I_r=I_lR_{r-l}$ for all $r\gg 0$. Then, for all $k\gg
0$, $X^k\in I$ and this implies that $l$ is even. On the other hand
$XY^kZ\in I$ for all $k\gg 0$; thus there exists $k_0$ such that
$XY^{k_0}Z\in I_l$ and $l=2+2k_0+3$, which means that $l$ is odd.
\end{example}

\subsection{Lexifiable ideals}


Let us consider a standard graded polynomial ring
$K[X_1,\ldots,X_n]$ with the degree lexicographic order. 
We recall that a {\it lexsegment (of degree $d$)} is a set $L$ of
monomials of degree $d$ with the property: if $u\in L$ and $v>u$ with
$\deg v=\deg u$ then $v\in L$. A homogeneous monomial ideal is said to be lexicographic if
all its graded components are spanned as a $K$-vector space by
lexsegments. It is clear that these definitions can be overtaken and
used in the weighted case.

\noindent
Lexicographic ideals play a central role in many 
results in commutative algebra because of their well understood
structure.  It would be of great interest to grasp which of
their properties also hold in the weighted case.
One of the most important facts concerning a lexsegment $L$ is that the
so-called {\it shadow}, i.e. the set of monomials which are obtained by
multiplying the monomials of $L$ by all the variables, 
is still a lexsegment. In the weighted case, given a lexsegment $L$ it
is natural to consider the $n$ shadows $L{\bf X_1},\ldots,L{\bf X_n}$.
In general they are not lexsegments, as the following easy example shows.


\begin{example}\label{ehbeh}
Let $(R,w)=(K[X,Y],(2,3))$, and consider the
monomial $XY$. This is the only monomial of degree $5$, and $\{XY\}$ is
obviously a lexsegment. Its shadow in degree $8$ is $\{XY^2\}$, which
is not lexsegment, since $X^4$ does not belong to it. 
\end{example}

It is thus quite clear that there are strong restrictions also on Hilbert
functions of lexsegment ideals generated in one degree. For instance,
an ideal $I$ generated in degree $d$ is a lexicographic ideal only if
contains $X_{11}^k$ 
for some $k>d$; this is possible only if $d=\alpha q_1$.

\begin{example}
Let $I\subseteq R=K[{\bf X_1},\ldots,{\bf X_n}]$ with $n>1$  be an
ideal generated by a lexsegment in degree $d=\alpha q_1$ and $H_I(d)\leq
l_1$. Then $I$ is a lexicographic ideal.\\
In fact $I$ is  generated by the
lexsegment $\{X_{11}^d,X_{11}^{d-1}X_{12},\ldots,X_{11}^{d-1}X_{1h}\}$, with
$h\leq l_1$; if we let $m\in I_r$ be a monomial with $r\geq
d$ then, since $m$ belongs to $I$, $m=X_{11}^{d-1}X_{1j}m\p$ for some
$j=1,\ldots,h$ 
and a monomial $m\p\in R_{r-d}$. If $s$ is  a monomial with $\deg s=\deg m$ and
$s\geq m$, then $s$ must be $X_{11}^{d-1}X_{1j\p}s\p$ with $j\p<j$ or
$j\p=j$ and $s\p\geq m\p$. It is thus clear that $s\in I_r$.
\end{example}

\begin{example}
Let $(R,w)=(K[X,Y,Z],(2,2,3))$. By the previous example, the ideal $(X^3,X^2Y)$ is
lexicographic, whereas $(X^3,X^2Y,XY^2)$ is not. However the ideal
$(X^3,X^2Y,XY^2,X^2Z^2)$ is lexicographic, as an easy verification shows.
\end{example}

At this point it is still not evident if it is possible to determine
whether an ideal is lexicographic  by looking at finitely many of its
graded components. The following proposition yields that this is
indeed the case. 

\begin{proposition}\label{chebene}
Let $I\subset (R,w)$ be a homogeneous ideal generated
in degree $\leq d$ and let $q\doteq\lcm(q_1,\ldots,q_n).$ If $I_i$ is
spanned
(as a $K$-vector space) by a lexsegment for all $i\leq d+q+G(w)$, then
$I$ is a lexicographic ideal.
\end{proposition}
\begin{proof}
The proof is by induction on $i$. We only have to prove that $I_i$ is
spanned
by a lexsegment if $i> d+q+G(w)$, provided that this is true for $I_r$
with $r<i$. If $I_i=\emptyset$ there is nothing to prove. Else, let
$v_i$ be the smallest monomial in $I_i$ and let $u>v_i$, with
$\deg u =\deg v_i$, be a monomial not in $I$. Finally let $X_{jh}$
denote the (lex-)smallest variable which divides $u$. Now we write
$v_i$ as $vm$, where $v$ is a minimal generator of $I$, we
say of degree $d\p \leq d$, and $m$ is the smallest monomial in
$R_{i-d\p}$.
Since $i-d\p>q+G(w)$, by Proposition \ref{beltrametti}, we may write
$m=m\p m\ppr$, where $m\p$ is the smallest monomial of $R_q$, which is
$X_{nl_n}^{q/q_n}$. Thus, $v_i=vX_{nl_n}^{q/q_n}m\ppr$. If we now let
$w\doteq
\frac{v_iX_{jh}^{q/q_j}}{X_{nl_n}^{q/q_n}}$, it is clear that $u\geq w\geq v_i$,
$w\in I_i$ and $w/X_{jh}\in I_{r-q_j}$. But this is a contradiction,
since $u/X_{jh}\geq  w/X_{jh}$ and $I_{r-q_j}$ is spanned by a lexsegment.
\end{proof}

In a polynomial ring with two variables and coprime weights, one
can expect to have a description of lexicographic ideals, because of
the following observation. Given a polynomial ring $(R,w)$  
according to \cite{CL} we call any set of consecutive monomials of the
same
 degree {\it a block}. If $R$ is a polynomial ring in two variables, any shadow of
a block is clearly a block. With this notation, a lexsegment is an
initial block, and Example \ref{ehbeh} shows that in general the
shadow of an initial block needs not to be such.

\noindent
Before proceeding with the characterization of lexicographic
ideals  of $K[X,Y]$,   we need to fix some notation. 
Given any set $A$  in $R$ of monomials of degree $d$, we let 
$\Shad_i(A)\subseteq R_{d+i}$ denote the set of the elements $um$, where $u\in A$
and $m$ is a monomial in $R_i$. Clearly, the
cardinality of $\Shad_{q_1}(A)$ equals that of $A$. Moreover, 
$|\Shad_{q_2}(A)|=|A|$ if $q_1\neq 1$ and $|\Shad_{q_2}(A)|=|A|+1$ if
$q_1=1$ and $A$ is a block.

\noindent
Finally, if $d\in\NN$ and $q_1\nmid d$,
we let $\delta$ denote the smallest integer $d+\beta q_2$,
$\beta\in\NN$, divisible by $q_1$. It is not difficult to see
that such a number exists and it is such that $\beta <q_1$ since $q_1$ and $q_2$ 
are assumed to be coprime.

\begin{lemma}\label{tecnico}
Let $L $ be a lexsegment of degree $d$.
\begin{itemize}
\item[(i)]
If $q_1\mid d$ then $Shad_i(L)$ is a lexsegment for all $i$.
\item[(ii)]
If $q_1\nmid d$ then $\{X^{\delta/q_1}\}\sqcup \Shad_{\delta-d}(L)$ is a
lexsegment (of degree $\delta$).
\end{itemize}
\end{lemma}
\begin{proof}
The proof of $(i)$ is obvious.

\noindent
To prove $(ii)$, let $X^aY^b$ be the largest monomial of $L$ and $\delta=d+\beta q_2$.
First observe that $b<q_1$.
Secondly, notice that  the largest monomial of
$\Shad_{\delta-d}(L)$ is $X^aY^{b+\beta}$. Since $b+\beta<2q_1$ and is
a multiple of $q_1$, we have $b+\beta=q_1$, so that the only monomial
of $R_\delta$ which is larger is $X^{\delta/q_1}$.
\end{proof}

\begin{theorem}\label{pesiprimi}
Let $I$ be a monomial ideal minimally generated in degrees $d_1< d_2 <
\ldots < d_r$ such that the monomials of $I_{d_i}$ form a lexsegment
for all $i=1,\ldots,r$. Then $I$ is a lexicographic ideal if and only
if $q_1\mid d_1$ or $q_1\nmid d_1$ and there exists $1<s\leq r$
such that $q_1\mid d_s$ and $d_s \leq\min\{\delta_i\}$.
\end{theorem}
\begin{proof}
We start by proving that if $q_1\mid d_1$ then $I$ is a lexicographic
ideal. For this purpose it is enough to observe that the shadow of a
block is a block and use Lemma \ref{tecnico}.

\noindent 
We show now that if $q_1 \nmid d_1$, the conditions on
$d_s$ imply that $I$ is lexicographic. In fact it is sufficient to
verify that $\Shad_i(\{I_{d_j}\})$ are lexsegments for all $i$ and
$j=1,\ldots,r$. Since the generator $X^{d_s/q_1}$ occurs in degree
$\leq\min\{\delta_i\}$, the conclusion follows again by Lemma \ref{tecnico}.

\noindent
Finally, if $I$ is lexicographic and $q_1 \nmid d_1$, then 
$X^k\in I_{kq_1}$ with $kq_1=d_s$ for some $1<s\leq r$. By Lemma
\ref{tecnico} $(ii)$ it is thus clear that $d_s\leq \min\{\delta_i\}$.
\end{proof}

\noindent
The conditions in the previous proposition can be easily re-formulated
in terms of Hilbert series. In general, it would be interesting
to have a solution for the following problem.

\begin{problem}\label{problema}
Find a combinatorial characterization for the Hilbert series of 
lexicographic ideals.
\end{problem}

In the same spirit of \cite{MP}, we say that an ideal $I\subseteq
(R,w)$ is {\it
lexifiable} if there exists a lexicographic ideal $L$ with the same
Hilbert function as $I$.

\noindent
Given a subset $A$ of monomials in $R_i$, $i\in\NN$, we let $\LEX(A)$
denote the set of the $|A 
|$ lexicographic largest elements of $R_i$. We also let 
$L\doteq\oplus_{i\in\NN}\langle\LEX(\{I_i\})\rangle$. 
Thus, $I$ is lexifiable iff $L$ is an ideal of $(R,w)$.
To establish which ideals are lexifiable is not an easy task. The
following example shows an ideal which is not lexifiable in any
lex-order. 

\begin{example}\label{lexi}
Let $(R,w)=(K[X,Y],(2,3))$. The ideal $(Y)$ provides an easy example of an
ideal which is lexifiable if $Y>_{\LEX} X$ and not lexifiable if
$X>_{\LEX} Y$.\\
Let $I=(X^3Y^3,X^2Y^4)$. $I$ is not lexifiable
in both cases $X>_{\LEX} Y$ and $Y>_{\LEX} X$. If $X>_{\LEX} Y$ then the candidate to be the
associated lexicographic ideal with $I$ is the ideal $L=(X^8,X^6Y)$ but
$H_I(18)=1$ and $H_L(18)=2$. If $Y>_{\LEX} X$ the
candidate is $L=(Y^5,Y^4X^2)$, but again $H_L(18)=2$.
\end{example}

\begin{example}
Let  $(R,w)=(K[X,Y],(2,7))$. The monomials of degree $28$ and $35$ are
$X^{14},X^7Y^2,Y^4$ and $X^{14}Y,X^7Y^3,Y^5$ respectively. Let us consider
the ideals which have exactly one minimal generator in these two
degrees. These are $I_1=(X^{14},X^7Y^3)$, $I_2=(X^7Y^2,Y^5)$, $I_3=(X^{14},Y^5)$,  
$I_4=(X^7Y^2,X^{14}Y)$, $I_5=(Y^4,X^7Y^3)$ and
$I_6=(Y^4,X^{14}Y)$. According to our definitions $I_1$ is
a lexicographic ideal, $I_2$ is lexifiable and $I_1$ is the lexicographic
ideal associated with it, $I_3$ is lexifiable associated with
$(X^{14}, X^7Y^3, Y^7)$. The ideals $I_4$, $I_5$ and $I_6$ are not
lexifiable, as a computation of the Hilbert function in degree $42$,
$42$ and $56$ shows.
\end{example}

Again according to \cite{MP}, we say that a graded polynomial ring $(R,w)$ is
{\it Macaulay-Lex} if every homogeneous ideal in $(R,w)$ is lexifiable.
Macaulay's Theorem together with Remark \ref{tuttouguale} says that 
$(R,w)$ with $w=(a,\ldots,a)$ is Macaulay-Lex, whereas for a general choice of
$w$ there are many ideals which are not lexifiable. Thus it is natural
to ask the following question.

\begin{question}\label{acaulay}
Which  polynomial rings $(R,w)$ are Macaulay-Lex?
\end{question}

\noindent
The results of the last part of this section shed some light on the problem
and provide partial answers to the above question.

\begin{theorem}\label{pesi1}
Let $I$ be a homogeneous ideal in $(R,w)=(K[X,Y],(1,q_2))$. There
exists a unique lexicographic ideal $L$ such that $H_{R/I}(t) = H_{R/L}(t)$ for
any $t \in \NN$. 
\end{theorem}
\begin{proof}
Taking in consideration what we have said before Example
\ref{lexi}, 
it is sufficient to prove that
$$\Shad_{1}(\LEX(\{I_d\}))\subseteq \LEX(\{I_{d+1}\}) \hbox{ and }  
\Shad_{q_2}(\LEX(\{I_d\}))\subseteq \LEX(\{I_{d+q_2}\}).$$ Since
$q_1=1$, $\Shad_i$ of an initial block is still an initial block,
and therefore we can reason on cardinalities.

\noindent
The first inclusion is immediate since, for any $A$, 
$|\Shad_1(A)|=|A|$, and consequently $|\Shad_1(\LEX(\{I_d\}))| =
|\LEX(\{I_d\})| = |\{I_d\}| = |\Shad_1 (\{I_d\})|$ which is equal to
$|\LEX(\Shad_1 (\{I_d\}))|$.\\
For the second inclusion,  we write $\{I_d\}$ as $\sqcup_{i=1}^s B_i$,
where $B_i$ are  maximal blocks. It is easy to
see that $\Shad_{q_2}(B_i) \cap \Shad_{q_2}(B_j) = \emptyset$.
Therefore
\begin{equation*}\begin{split}
|\LEX(\Shad_{q_2} (\{I_d\}))| - |\LEX(\{I_d\})|& = |\Shad_{q_2} (\{I_d\})| -
|\{I_d\}|\\
& =\sum_{i=1}^s |\Shad_{q_2}(B_i)| - |B_i|\\
&\geq 1\\
& =|\Shad_{q_2}(\LEX(\{I_d\}))| - |\LEX(\{I_d\})|,
                 \end{split}
\end{equation*}
and the proof is concluded.
\end{proof}

\begin{example}
This is an example of an ideal $I$ which is not lexifiable in a ring
for which Condition \ref{multipli} is satisfied.\\
Consider $(R,w)=(K[X,Y,Z],(1,2,4))$ and let $I=(X^4,Y^2,X^3Y)$, for which we have
that $H_I(4)=2$, $H_I(5)=3$, $H_I(6)=4$, $H_I(7)=4$. $I$ is not
lexifiable, in fact, if we try to
construct the associated lexicographic ideal $L$, we shall have 
$L_4=\{X^4,X^2Y\}$, $L_5=\{X^5, X^3Y, XY^2\}$ and $L_6=\{X^6, X^4Y, X^2Y^2,
X^2Z\}$, so that $H_L(7)\geq 5$.
\end{example}

\begin{example}
Let $(R,w)=(K[X,Y,Z],(1,a,ab))$. The ideal
$$I=(X^{ab},Y^b,X^{a+1}Y^{b-1})$$ is not lexifiable, therefore $(R,w)$
is not Macaulay-Lex.
\noindent
Let $b>2$ and suppose that $I$ is lexifiable with associated lexicographic ideal
$L$. 
We first observe that  $I_j$ does not contain any monomial divisible by $Z$
for all $j<2ab$. Secondly, we show that $X^{ab-2a}Z$ is a monomial of $L$.
Since $H_I(ab+(\alpha+1) a)\geq H_I(ab+\alpha a)+1$ for all
positive $\alpha\in\ZZ$, and $H(ab+a)=5$, one gets that
$H_I(ab+(b-2)a)\geq
H(ab+a)+b-3=b+2$. Since the first $b+2$ monomials in degree $ab+\alpha
a$ are $X^{ab+\alpha a},X^{ab+(\alpha-1)a}Y,\ldots,X^{\alpha
  a}Y^b,X^{\alpha a}Z$, this proves our claim.  
It is  convenient now to write down all the monomials of degree $2ab-2a$
which are $\geq Y^{2b-2}$, which is the smallest monomial of $I$ in
this degree. These are
\begin{equation*}\begin{split}
&X^{2ab-2a},X^{2ab-3a}Y,\ldots,X^{ab-2a}Y^b,X^{ab-2a}Z,\\
&X^{ab-3a}Y^{b+1},X^{ab-3a}YZ,\ldots,X^aY^{2b-3},X^aY^{b-3}Z,
Y^{2b-2}.        \end{split}
\end{equation*}

\noindent
As a consequence, it is easy to compute that 
$H_I(2ab-2a)$ is 
$(b+1)+(b-3+1)=2b-1$. Analogously one gets that $H_I(2ab-1)=2b$. Furthermore
the monomials of $L_{2ab-2a}$ have a certain number of multiples in
degree $2ab-1$, which we can count. We can thus estimate the
cardinality of 
$L_{2ab-1}$ as follows:
\begin{equation*}\begin{split}
|L_{2ab-1}|&\geq
|\{u\in L_{2ab-1} \: Z\nmid u\}
+|\{u\in L_{2ab-1} \: Z\mid u\}|\\
&=\left[(b+1)+\left\lceil\frac{b-3}{2}\right\rceil+1\right]
+\left[\left(\left\lfloor\frac{b-3}{2}\right\rfloor+1\right)+1\right]\\
&=2b+1.
                 \end{split}
\end{equation*}
This is a contradiction since $I$ and $L$ have by definition the same
Hilbert function.

\noindent
Finally, if $b=2$, an easy computation of the Hilbert function in
degree $3a+1$ shows that $I$ is not lexifiable.
\end{example}

\subsection{Polarization}

In \cite{P} it is  shown how in the standard case the 
lexicographic ideal $L$ associated with an ideal $I\subseteq R$ can be
also obtained as  the result of a finite process which consists of
three 
fundamental steps, which are a) polarizing a monomial ideal  b) modding out by a
sequence of generic linear forms and c) taking initial ideals (with
respect to the lexicographic order). In the non-standard case,
following step-by-step the original proof and  using generic
sequences of homogeneous forms (which are not necessarily linear) 
it is not difficult to prove that the same procedure also
terminates and leads to an ideal $I\pol$, which we call {\it completely polarized} . Since, as
we have already seen, not all ideals are lexifiable, one might make
use of the $I\pol$, which is a strongly stable monomial ideal with the
same Hilbert function as $I$, instead.

\begin{example}
Let $(R,w)=(K[X,Y,Z],(1,2,4))$  and  
$$I=(X^8,X^6Y,X^4Y^2,X^2Y^3,Y^4,X^2YZ,X^6Z).$$
One can verify that the ideal $$L=(X^8,X^6Y,X^4Y^2,X^4Z,X^2Y^3,X^2YZ,X^2Z^2,Y^6)$$ is 
the lexicographic ideal associated with $I$ and,
thus, $I$ is lexifiable. On the other hand 
$$I\pol=(X^8,X^6Y,X^4Y^2,X^4Z,X^2Y^3,X^2Y^2Z,Y^4).$$
This shows that, even in the favourable case when Condition \ref{multipli} is
satisfied, one might have  $I\pol\neq L$.
\end{example}

We conclude by posing the following question.

\begin{question}
Is there a combinatorial characterization  of completely polarized ideals?
\end{question}


{\small
\begin{thebibliography}{11111111}

\bibitem [A] {A} V. Arnold. ``$A$-graded algebras and continued
  fractions''. {\em Communications in Pure and Applied Math.} {\bf
  42} (1989), 993-1000.
\bibitem [ADK] {ADK} A. Aramova, E. De Negri, K. Crona. ``Bigeneric
  initial ideals, diagonal subalgebras and bigraded Hilbert
  functions''. {\em J. Pure Appl. Algebra} {\bf 150} (2000), No. 3, 215-235.
\bibitem [B] {B} V. Bavula. ``Identification of the Hilbert function
  and Poincare Series and the dimension of modules over filtered
  rings''. {\em Russian Acad. Sci. Izv. Math.} {\bf 44} (1995), No. 2, 225-246.
\bibitem [Be] {Be} D. Benson. ``Dickson invariants, regularity and computation in group
cohomology''. 
{\em Illinois J. Math.} {\bf 48} (2004), No. 1, 171-197.
\bibitem [BH] {BH} W. Bruns and J. Herzog. 
{\em Cohen-Macaulay rings}. Cambridge University Press, Cambridge,
1998.
\bibitem [BR] {BR} M. Beltrametti and L. Robbiano. 
``Introduction to   the theory of weighted projective spaces''. {\em Expo. Math.} {\bf 4}
  (1986), 111-162.
\bibitem [BS] {BS} D. Bayer and M. Stillman. 
``A theorem on refining division orders by the
reverse lexicographic orders.'' {\em Duke J. Math.} {\bf 55} (1987), 321-328.
\bibitem [CL]  {CL} G. Clements and B. Lindstr\"om. ``A
  generalization of a combinatorial theorem of Macaulay''. {\em
  J. Combinatorial Theory} {\bf 7} (1969), 230-238.
\bibitem [CoCoA] {CoCoA} CoCoA: a system for doing 
     Computations in Commutative Algebra,
  Available at http://cocoa.dima.unige.it.
\bibitem [D] {D} C. Delorme. ``Espaces projectifs anisotropes''. {\em
  Bull. Soc. Math. France} {\bf 103} (1975), 203-223.
\bibitem [DH] {DH} E. De Negri and  J. Herzog. ``Completely lexsegment ideals''.
{\em Proc. Amer. Math. Soc.} {\bf  126}  (1998)  No. 12, 3467-3473.
\bibitem [E] {E} D. Eisenbud. 
{\em Commutative Algebra}. Springer-Verlag, New York, 1995.
\bibitem [EK] {EK} S. Eliahou and M. Kervaire, Minimal resolutions of
  some monomial ideals, J. Algebra {\bf 129} (1990), 1-25.
 \bibitem[M] {M}
F. Macaulay. Some properties of enumeration in the theory of modular systems.
{\em Proc. London Math. Soc.} {\bf 26} (1927), 531-555.
\bibitem [MP] {MP} 
J. Mermin and I. Peeva. ``Lexifying ideals''. Preprint 2005.
\bibitem [P] {P}
K. Pardue. {\em Nonstandard Borel-fixed ideals}. Thesis, Brandeis
University, 1994.
\bibitem [S] {S} B. Sturmfels. {\em Gr\"obner Bases and Convex
  Polytopes}. University Lecture Series, Vol. 8, Amer. Math. Soc.,
  Providence, 1996.
\bibitem [SS] {SS} G. Scheja and U. Storch. {\em Regular Sequences and
  Resultants}. Research notes in Mathematics, Vol. 8, A K Peters,
  Natick, Massachusetts, 2001.
\end{thebibliography}}
\end{document}